\documentclass{article}

\usepackage{latexsym}
\usepackage[centertags]{amsmath}
\usepackage{amsfonts}
\usepackage{amssymb}
\usepackage{amsthm}
\usepackage{newlfont}
\newcommand{\Real}{\mathbb R}
\newcommand{\Complex}{\mathbb C}
\newcommand{\Proj}{\mathbb P}

\begin{document}

\begin{center}
{\LARGE \bf On the algebraic invariant curves of plane polynomial
differential systems}\\

 \vspace{0.3cm}

\ {\Large \bf  Tsygvintsev Alexei}

\end{center}

\begin{abstract} We consider a plane polynomial vector field $P(x,y)dx+Q(x,y)dy$ of degree
$m>1$. To  each algebraic invariant curve of such a field we
associate a compact Riemann surface with the meromorphic
differential $\omega=dx/P=dy/Q$. The asymptotic estimate of the
degree of an arbitrary algebraic invariant curve is found. In the
smooth case this estimate was already found by D. Cerveau and A.
Lins Neto [2] in a  different way.
\end{abstract}

\begin{center}
{\bf Introduction}
\end{center}

The study of plane polynomial vector fields goes back at least to
Poincar\'e [12]. Recall that the second half of Hilbert's 16th
problem [8] asks for an upper bound on the number of limit cycles
of real plane polynomial vector fields. Notice, that the class of
invariant curves of the given planar system  involves the class of
its limit cycles. Of course, every limit cycle is also an
invariant curve.

This paper is devoted to one aspect of this problem: to study
algebraic invariant curves i.e. defined by an  algebraic equation
$f(x,y)=0$, where $f\in \Complex [x,y]$ is an arbitrary
polynomial. The real part of the above curve which turns out to be
a limit cycle, is called the {\it algebraic} limit cycle. Up to
now only several cases of algebraic limit cycles
 are known, especially for quadratic plane systems [3].
It has been shown by Darboux [4] that if a given planar polynomial
system of degree $m$ has more than $2+[m(m+1)]/2$ algebraic
invariant curves, then it admits a rational first integral.

In this paper we apply a new method connecting the problem of
existence of algebraic invariant curves of plane polynomial vector
fields of the form $P(x,y)dx+Q(x,y)dy$ with the contemporary
theory of Riemann surfaces. To each algebraic invariant curve of
such a field we associate a compact Riemann surface $C$ and a
meromorphic differential $\omega=dx/P=dy/Q$.

Using this approach, in Section 5 we  find  the asymptotic
estimate of the degree of an arbitrary algebraic invariant curve
(Theorem 6). In the particular case we obtain the estimate for a
degree of a nodal algebraic invariant curve (Corollary 3). It is
shown too that an arbitrary smooth algebraic invariant curve has a
degree less than $m+2$ (Theorem 2) and that for an arbitrary
algebraic invariant curve its genus is a linear function of the
degree (Theorem 5). These results were already obtained (in a
completely different way) in papers [1], [2].

\begin{center}
{\bf 1. The Darboux divisor and points at infinity.}
\end{center}

Consider the system of differential equations $$ \dot x=P(x,y),
\quad \dot y=Q(x,y), \quad (x,y)\in \Complex^2 \eqno (1) $$ where
$P,Q$ are polynomials of degree $m>1$. We suppose that $P$ and $Q$
have not a common nonconstant polynomial factor and
$P=\sum\limits_{i=1}^mP_i$, $Q=\sum\limits_{i=1}^mQ_i$, where
$P_i$,$Q_i$ are homogeneous polynomials of degrees $i=0,\dots,m$.
\vspace{0.3cm}

Let $C=\{(x,y)\in \Complex^2 :f(x,y)=0\}$ be an invariant curve of
(1). Without loss of generality we may suppose that $f\in \Complex
[x,y]$ is irreducible. Then $\dot
f=\left(P\displaystyle\frac{\partial f}{\partial
x}+Q\displaystyle\frac{\partial f}{\partial y}\right)_{f=0}\equiv
0$. As the ideal $<f>$ is radical, then $\dot f\in <f>$ and hence
$\dot f=kf$, for some $k\in  \Complex [x,y]$. \vspace{0.3cm}

\noindent{\bf Definition 1}  The polynomial $f(x,y)\in \Complex
[x,y]$ is called an algebraic partial integral of the system (1)
if there exists a polynomial $k\in \Complex [x,y]$ such that $$
P\frac{\partial f}{\partial x}+Q\frac{\partial f}{\partial y}=kf.
\eqno (2) $$ The polynomial $k$ is called {\it cofactor} and has
the following form $k=\sum\limits_{i=1}^{m-1}k_i$, where $k_i$ are
homogeneous polynomials of degrees $i=0,\dots,m-1.$\\

If $k\equiv 0$ then $f(x,y)=\mathrm{const}$ is a first integral of
the system (1).

\vspace{0.3cm}

\noindent{\bf Remark 1} It is easy to see that if $f(x,y)$ is
reducible,  i.e. $f=f_1^{m_1}\cdots f_l^{m_l}$, where $f_k\in
\Complex [x,y]$, $k=1,\dots,l$, then polynomials $f_k$ are again
partial integrals of the system (1).

The polynomial $f$ is a sum of its homogeneous parts
$f=\sum\limits_{i=a}^nf_i$, where $f_i$ are homogeneous
polynomials of degrees $i=0,\dots,n$ and  $n=\mathrm{deg}(f)$.

Consider the homogeneous polynomial $R_{m+1}(x,y)$ of degree $m+1$
defined by $$ R_{m+1}(x,y)=xQ_m(x,y)-yP_m(x,y), \eqno (3) $$ where
$P_m$ and $Q_m$ are higher homogeneous parts of the polynomials
$P$ and $Q$ respectivelly. Let us suppose that $R_{m+1}$ does not
vanish identically, then it has $m+1$ zeros $D_i=[x_i:y_i]\in
{\Complex \Proj^1}$, $i=1,\dots,m+1$. By the suitable rotation of
variables $x$, $y$ we can obtain $x_iy_i\neq 0$, $i=1,\dots,m+1$.
Hence, without loss of generality: $D_i=(1,z_i)$, $z_i \in
\Complex^2$, $z_i \neq 0$, $i=1,\dots,m+1$. \vspace{0.3cm}

\noindent {\bf Definition 2} The formal sum of points $D=\sum\limits_{i=1}^{m+1}D_i$ is called the
{\it Darboux divisor} of the differential system (1).
\vspace{0.3cm}

Notice that the impotant role of the points $D_i$ for polynomial
vector fields first was observed by Darboux in 1878.

Let $$ V(x,y)=P(x,y)\frac{\partial}{\partial
x}+Q(x,y)\frac{\partial}{\partial y}, $$ be the polynomial vector
field on $\Complex^2$ corresponding to the system (1). Through the
non-linear change of variables $$ u=\frac{1}{x},\quad
v=\frac{y}{x}, \quad x\neq 0, \quad (u,v)\in  \Complex^2, $$ and
multiplying the induced vector field by $u^{m-1}$ we obtain [6],
[5] $$
\begin{array}{lll}
\tilde V(u,v)= A(u,v)\displaystyle\frac{\partial}{\partial
u}+B(u,v)\displaystyle\frac{\partial}{\partial v},\\
A(u,v)=-u^{m+1}P\left(\displaystyle\frac{1}{u},\displaystyle\frac{v}{u}\right),\\
B(u,v)=u^m\left[Q\left(\displaystyle\frac{1}{u},\displaystyle\frac{v}{u}\right)-
vP\left(\displaystyle\frac{1}{u},
\displaystyle\frac{v}{u}\right)\right],\\
\end{array}
$$
where $\tilde V(u,v)$ represents the vector field of (1) near the line at infinity
$L_\infty=\{u=0\}$. The point $(0,v_0)$ where $\tilde V(0,v_0)=(0,0)$  is the singular point of
$\tilde V(u,v)$. It is easy to see that $R_{m+1}(1,v_0)=0$ and we obtain
\vspace{0.3cm}

\noindent {\bf Proposition 1} The points $D_i=(1,z_i)\in D$, $i=1,\dots,m+1$ are the singular
points at infinity of the system (1).
\vspace{0.3cm}

The equation (2) turns into
$$
A(u,v)\displaystyle\frac{\partial F}{\partial u}+B(u,v)\displaystyle\frac{\partial F}{\partial
v}=K(u,v)F,
$$
where
$F(u,v)=u^nf\left(\displaystyle\frac{1}{u},
\displaystyle\frac{v}{u}\right)=f_n(1,v)+uf_{n-1}(1,v)+\cdots=0$  represents the curve $C$ near
$L_\infty$ and
$$
K(u,v)=u^{m-1}k\left(\displaystyle\frac{1}{u},
\displaystyle\frac{v}{u}\right)-u^mnP\left(\displaystyle\frac{1}{u},
\displaystyle\frac{v}{u}\right).
$$

 Let us show now that the Darboux divisor $D$ contains all possible  points at infinity of any
algebraic invariant curve of the system (1).

Denote by $L_\infty=\{[x_i:y_i:0]:(x,y)\subset \Complex \Proj
^1\}\subset \Complex \Proj^2$ the line at infinity. Let $I_f$ be a
set of points at infinity of the algebraic curve $C$ which
correponds to the
 equation $f(x,y)=0$, where $f(x,y)$ is an algebraic partial integral of the system (1).
\vspace{0.3cm}

\noindent{\bf Theorem 1.} $I_f\subset D$.
\vspace{0.3cm}

\noindent {\bf Proof.} By considering the right and left hand homogeneous parts of (2) we find
$$
P_m\frac{\partial f_n}{\partial x}+Q_m\frac{\partial
f_n}{\partial y}=k_{m-1}f_n, \eqno (4)
$$
where $f_n$ is the highest order term of the polynomial $f=\sum\limits_{i=a}^nf_i$ and $k_{m-1}$
is the highest order term of the
cofactor $k=\sum\limits_{i=1}^{m-1}k_i$.

To show $I_f\subset D$ we need to prove that if $f_n(x_0,y_0)=0$ then
$(x_0,y_0)\in D$ or
$$
R_{m+1}(x_0,y_0)=0, \eqno (5)
$$
where the polynomial $R_{m+1}$ is defined by (3).

 Consider the linear change of variables $(x,y)\rightarrow  (u,v)$: $ x=x_0+u$, $y=y_0+v$.
The polynomial $f_n(x,y)$ turns into the polynomial
$F(u,v)=f_n(x_0+u,y_0+v)$ which has the following Taylor expansion
$$ F(u,v)=\sum_{i=r}^nF_i(u,v),\eqno (6) $$ where $F_i$ are
homogeneous polynomials of degrees $i=r,\dots,n$, $r\geq 1$ and $$
F_i=\displaystyle\frac{1}{(n-i)!}\left(x_0\displaystyle\frac{\partial}{\partial
u}+y_0\displaystyle\frac{\partial}{\partial
v}\right)^{n-i}f_n(u,v). $$ Thus, for the  lower order term $F_r$
of the sum (6) we have $F_r\not\equiv \mathrm{const}$ and the
following identity is fulfilled $$ x_0\frac{\partial F_r}{\partial
u}+y_0\frac{\partial F_r}{\partial v}=0. \eqno (7) $$

The equation (4) takes the form
$$
((c_1+N_1(u,v))\frac{\partial}{\partial u}+(c_2+N_2(u,v))\frac{\partial}{\partial
v}\quad)(F_r+\cdots+F_n)=0, \eqno (8)
$$

where
$$
c_1=P_m(x_0,y_0)-\displaystyle\frac{x_0}{n}k_{m-1}(x_0,y_0),\quad
c_2=Q_m(x_0,y_0)-\displaystyle\frac{y_0}{n}k_{m-1}(x_0,y_0),
 \eqno (9)
$$
are constants and $N_1(u,v)$, $N_2(u,v)$ are polynomials such that $N_1(0,0)=N_2(0,0)=0$.

The two cases should be considered.

\noindent 1) $c_1=c_2=0$. Then from relations (9) it follows that the equality (5) is fulfilled.

Hence $(x_0,y_0,0)\in D$.

\noindent 2)  $(c_1,c_2)\neq(0,0)$. Then one can show from (8)
that $c_1\displaystyle\frac{\partial F_r}{\partial
u}+c_2\displaystyle\frac{\partial F_r}{\partial v}=0$. Using (7)
we see that vectors $(c_1,c_2)$ and $(x_0,y_0)$ are colinear i.e.
$$ \mathrm{det} \left(
\begin{array}{ll}
c_1&x_0\\
c_2&y_0\\
\end{array}
\right)
=0,
$$
which gives again the equality (5). Q.E.D.
\vspace{0.3cm}

\noindent {\bf Corollary 1.} Let $D=D_1,\dots,D_{m+1}$ be a
Darboux divisor of the system (1) and $l_i=a_ix+b_iy,$
$a_i$,$b_i\in \Complex$, $i=1,\dots,m+1$ be a set of linear forms
such that $l_i(D_i)=0$, $i=1,\dots,m+1$. Then there exists
nonnegative integers $n_1,\dots,n_{m+1}$,
 $\sum n_i=n$ that
$$ f_n(x,y)=\prod_{i=1}^{m+1}l_i^{n_i}(x,y).\eqno (10) $$ Notice,
that the same expression for $f_n$ was introduced first by
Jablonskii [9] in the case $m=2$, see also [10].

\begin{center}
{\bf 2. The smooth case.}
\end{center}

Let $C\subset \Complex \Proj^2$ be an algebraic smooth curve of
$\mathrm{deg}(C)=n$ satysfying the equation $f(x,y)=0$ where
$f(x,y)$ is an irreducible algebraic partial integral of the
system (1). Without loss of generality we suppose that $$
f(x,y)=y^n+a_1(x)y^{n-1}+\cdots+a_n(x), \quad a_i(x)\in \Complex
[x],\quad i=1,\dots,n. $$ Consider the holomorphic mapping
$\phi:C\rightarrow \Complex \Proj^1$ defined by $\phi(x,y)=x$.

Let $\nu=\nu_{\phi}(P)$ be a multiplicity of $\phi$ at the point
$P\in C$. Consider the ramification divisor $R=\sum\limits_{p\in
C}(\nu_{\phi}(P)-1)P\subset \mathrm{Div}(C)$.

We break $R$ into two divisors $R=R_1+R_2$, where $$
R_1=\sum\limits_{P\in C\cap L_\infty}(\nu_{\phi}(P)-1)P $$
contains branching points of $\phi$ at infinity and $$
R_2=\sum\limits_{P\in C/L_{\infty}}(\nu_{\phi}(P)-1)P $$ contains
all finite branching points. \vspace{0.3cm}

\noindent {\bf Lemma 1.} Let $C=\{f(x,y)=0\}\subset \Complex
\Proj^2$ be a nonsingular algebraic curve of $\mathrm{deg}(C)=n$
where $f(x,y)$ is a partial first integral of the system (1). Then
$$ \mathrm{deg}(R_1)\leq n-1. $$ \vspace{0.3cm}

This statement is proved by noting that $f=f_n+\cdots+f_0$,
$\mathrm{deg}f_k=k$ and $f_n=\prod\limits_{i=1}^m L_i^{n_i}(x,y)$
where $\sum\limits_{i=1}^m n_i=n$, $m\leq n$, $L_i(x,y)$ are
linear homogeneous polynomials. \vspace{0.3cm}

\noindent {\bf Lemma 2.}
$\mathrm{deg}(R_2)=n^2-n+1-\mathrm{deg}(R_1)$. \hspace{8.4cm} (11)
\vspace{0.3cm}

\noindent {\bf Proof.} Denote $g=\mathrm{genus}(C)$,
$n=\mathrm{deg}(C)$, then by the genus-degree formula for a
nonsingular curve $C$ we have
$g=\displaystyle\frac{(n-1)(n-2)}{2}$.\\ By the Riemann-Hurwitz
formula we obtain $g=\displaystyle\frac{\mathrm{deg}(R)}{2}-n+1$.
Comparing these two expressions for $g$ we find (11). Q.E.D.

Now let us study the divisor $R_2$.

If $K=(x_0,y_0)\in R_2$ then $\displaystyle\frac{\partial
f}{\partial y}(K)=0$ by the definition  of a branching point. With
help of (2) we obtain $$ P(K)\displaystyle\frac{\partial
f}{\partial x}(K)+Q(K)\displaystyle\frac{\partial f}{\partial
y}(K)=0. \eqno (12) $$ \vspace{0.3cm}

\noindent {\bf Lemma 3.} If the curve $C$ is nonsingular,
$\mathrm{deg}(C)=n$, then $$ \mathrm{deg}(R_2)\leq mn, $$ where
$m>1$ is the degree of the system (1). \vspace{0.3cm}

\noindent {\bf Proof.} Since $K$ is a smooth point the relation
(12) holds $$ \left(f,\frac{\partial f}{\partial y}\right)_K\leq
(f,P)_K,\quad K\in R_2, $$ where $(g,l)_X$ denotes the
intersection number of the curves $g(x,y)=0$ and $l(x,y)=0$ at the
point $X\in g\cap l$. One can easily verify that
$\mathrm{deg}(R_2)=\sum\limits_{P\in
R_2}\left(f,\displaystyle\frac{\partial f}{\partial y}\right)_P$.
Thus, by B\'ezout theorem $\mathrm{deg}(R_2)\leq mn$. Q.E.D.
\vspace{0.3cm}

\noindent {\bf Theorem 2.} Let us assume that the system (1)
admits an smooth algebraic invariant curve $C\subset \Complex
\Proj ^2$ defined by the equation $f(x,y)=0$, $\mathrm{deg}(f)=n$.
Then $$ n\leq m+1 $$ where $m>1$ is the degree of the  system (1).
\vspace{0.3cm}

The statement of the theorem follows immediately from the above
three lemmas. The Theorem 2 was obtained for the first time in [2]
using a different method. By J. Moulin-Ollagnier it was shown that
the same result can be obtained in the theory of the  Koszul
complexes of polynomial vector fields.

\begin{center}
{\bf 3. The Weierstrass polynomials.}
\end{center}

Let $X=(x_0,y_0)\in \Complex^2$ be a finite singular point of the
curve $C=\{(x,y)\in \Complex^2: f(x,y)=0\}$ i.e. the point $X$ at
which $\displaystyle \frac{\partial f}{\partial
x}(X)=\displaystyle \frac{\partial f}{\partial y}(X)=0$. Without
loss of generality we suppose $X=(0,0)$.

In order to clarify the local structure of $C$ near $X$, we shall
need the help of the {\it Weierstrass polynomials} [7n7].

Let $\Complex  \{x\}$ $\left( \Complex  \{x,y\} \right)$ represent
the ring of holomorphic functions defined in some neigborhood of
$0\in \Complex$ $((0,0)\in \Complex^2)$. \vspace{0.3cm}

\noindent {\bf Definition 3.} $w\in \Complex \{x,y\}$ is said to
be a Weierstrass polynomial with respect to $y$, if $$
w=y^d+c_1(x)y_{d-1}+\cdots+c_d(x),\quad c_j(x)\in \Complex  \{x\},
\quad c_j(0)=0, \quad j=1,\dots,d. $$ \vspace{0.3cm} Let us assume
that $C$ is irreducible and its affine equation is $$
f(x,y)=y^n+a_1(x)y^{n-1}+\cdots+a_n(x)=0. $$ \vspace{0.3cm}

\noindent {\bf Theorem 3.}  The polynomial $f(x,y)$ can be
expressed as $$ f=uf_1 f_2 \cdots f_p, $$ where
$f_i(x,y)=y^{d_i}+c_{i1}(x)y^{d_i-1}+\cdots+c_{id_i}(x)$,
$i=1,\dots,p,$ are irreducible Weierstrass polynomials and
$u(x,y)$  is a unit of $\Complex \{x,y\}$, i.e. $u(0,0)\neq 0$.

There exists the open discs $\Delta_i=\{\tau\in \Complex: \mid
\tau \mid<\rho_i\}$, $i=1,\dots,p$, such that each equation
$f_i(x,y)=0$, $i=1,\dots,p$ defines holomorphic mapping
$q_i:\Delta_i \rightarrow C$ as follows $$ \tau\rightarrow \left(
\tau^{d_i},g_i(\tau) \right), \quad   \mathrm{where} \quad
g_i(\tau)=\sum_{k=1}^\infty c_{ik} \tau^k \in \Complex \{\tau\},
\quad i=1,\dots,p. \eqno (13)$$ \vspace{0.3cm}

Thus, with topological point of view, the algebraic curve $C$ can
be obtained near the singular point $X=(0,0)$ from several open
discs by identifying them together at their centers. This is the
concept of normalization [7].

\vspace{0.3cm}

\noindent {\bf Theorem 4.} Let $C=\{(x,y)\in \Complex \Proj^2:
f(x,y)=0\} $ be an algebraic invariant curve of the system (1) and
$X=(x_0,y_0)$ be a singular point of $C$. Then $X$ is an
equilibrium point of the system (1). \vspace{0.3cm}

\noindent {\bf Proof.} Let us assume that $X=(x_0,y_0)$ is not an
equilibrium point point of the system (1). Then it has the unique
solution passing through this point $$
x=x_0+P(x_0,y_0)t+\sum\limits_{i=2}^\infty a_it^i,\quad
y=y_0+Q(x_0,y_0)t+\sum\limits_{i=2}^\infty b_it^i, \quad a_i,  b_i
\in \Complex \eqno (14) $$ where $t\in\Delta=\{t\in \Complex :\mid
t \mid < \rho\}$ for any small $\rho\in \Real$.

On the other hand $X$ is the singular point of $C$ and according to the Theorem 3 the
system (1) has no less than $p>0$ different solutions passing through $X$ and locally expressed by
(13).
Thus, we obtain $p=1$ and the solution (14) is the parametrization of the curve $C$ near the
singular point $X$. By our assumption $X$ is not an equilibrium point point of (1) i.e.
$P(x_0,y_0)\neq 0$ or $Q(x_0,y_0)\neq 0$.  Hence , looking at (14), $X$ is the smooth point of
$C$. We obtain the contradiction. Q.E.D.
\vspace{0.3cm}

\noindent {\bf Corollary 2.} The number of finite singular points
of an arbitrary algebraic invariant curve of the system (1) is not
more than $m^2$. Furthermore, if $\displaystyle
\frac{P_m(x,y)}{Q_m(x,y)}\not \equiv \displaystyle\frac{x}{y}$,
then $$ \mid \mathrm{Sing}(C)\mid\leq m^2+m+1. $$ \vspace{0.3cm}

Indeed, if  $\displaystyle\frac{P_m(x,y)}{Q_m(x,y)}\not \equiv
\displaystyle\frac{x}{y}$ then the polynomial (3) is not equal
zero identically and according to Corollary 1 the curve $C$ cannot
have more than $m+1$ singular points at infinity.

\begin{center}
{\bf 4. The genus of $C$.}
\end{center}

Let $C$ be an algebraic invariant curve of the system (1) defined
by the equation $f(x,y)=0$. Denote by $\mathrm{Sing}(C)$ the set
of its singular points. There exists the compact Riemann surface
$\tilde C$ with a surjective continuous map $\pi:\tilde C
\rightarrow C$ such that $\pi : \tilde C
/\pi^{-1}(\mathrm{Sing}(C))\rightarrow C/\mathrm{Sing}(C)$ is a
holomorphic bijection. The aim of this section is to calculate the
genus  of $\tilde C$ which is also called the genus of the curve
$C$. Consider the following meromorphic differential on $C$ $$
\omega=\displaystyle\frac{dx}{P}=\displaystyle\frac{dy}{Q}. \eqno
(15) $$

Let $\omega$ be its divisor then according to the Poincar\'e-Hopf
formula $$ 2g-2=\mathrm{deg}(\omega).  \eqno (16) $$ On the other
hand, by Noether's formula [8n11] $$
g=\displaystyle\frac{(n-1)(n-2)}{2}-\sum_{X\in \mathrm{Sing}(C)}
\delta (X), \eqno (17) $$ where the numbers $\delta(X)$ are given
by $$ \delta(X)=\left(f,\frac{\partial f}{\partial
y}\right)_X+\mid\pi^{-1}(X)\mid-\nu_{\phi}(X). $$ Here $( , )_X$
is the intersection number and $\nu_{\phi}(X)$ is the multiplicity
of the map $\phi: (x,y)\rightarrow x$ at the point $(x,y)\in
\mathrm{Sing}(C)$.

It is easy to see that $\omega$ has no zeros in the affine part of
$C$. Let now $ X=(x_0,y_0)\in \Complex^2$ be the singular point of
the curve $C$. Without loss of generality we put $X=(0,0)$.
According to Theorem 3 we can factor $f(x,y)$ into the product of
irreducible factors $$ f=uf_1\cdots f_r, $$ where $u(0,0)\neq 0$
and $f_i$, $i=1,\dots,r$ are Weierstrass polynomials. Notice, that
$\mid\pi^{-1}(X)\mid=r$. Then locally $C$ can be represented as
follows $$ C=C_1+\cdots+C_r, $$ where $C_i=\{(x,y)\in
\Complex^2:\mid x \mid < \rho, \mid y \mid < \epsilon,
f(x,y)=0\},\quad i=1,\dots,r $ are irreducible local analytic
curve components of $C$ and  $\rho$, $\epsilon$ are sufficiently
small real numbers.

The parametrization of $C_i$,  $i=1,\dots,r$  near $X=(0,0)$  is
given by $$ x=\tau^{d_i},\quad y=\sum^\infty_{k=1}c_{ik}\tau^k,
\quad c_{ik}\in \Complex ,\quad d_i=\mathrm{deg}(f_i).\eqno (18)
$$

Puting (18) into (15) and using Theorem 4  one can show that the
differential $\omega$ has in the point $X$ a pole of the
multiplicity at least one. So, for the affine part of the curve
$C$ we have the following estimate $$ \mathrm{deg}(\omega) \mid_{C
\cap \Complex^2 } \leq -\sum_{X\in \mathrm{Sing}(C)\cap \Complex
^2} \mid \pi^{-1}(X) \mid. \eqno (19) $$

Now let us consider the points at infinity. Substituting $x=1/u$,
$y=v/u$ into $f(x,y)=0$ and multiplying both sides of the
resulting expression by $u^n$, we obtain the equation $$
F(u,v)=f_n(1,v)+uf_{n-1}(1,v)+\cdots+f_0u^n=0, \quad
f_0=\mathrm{const}\neq 0 , $$ which represents the algebraic curve
curve $C$ near the line at infinity $L_\infty=\{u=0\}.$ We can
write $f_n(1,v)$ as follows $$
f_n(1,v)=\prod_{i=1}^q(v-v_i)^{n_i},\quad n_i=0,1,\dots,\quad
q\leq n,\quad \sum n_i=n,\eqno (20) $$ where the points
$(0,v_i)\in C\cap L_{\infty}$, $i=1,\dots,k$.

Now we break (20) into the product of three factors $$
f_n(1,v)=L_1L_2L_3. $$ Here $L_1=\prod\limits_{i=1}^r(v-v_{1i})$,
$r\leq n$ contains all simple factors of (17). Near the points
$(0,v_{1i})$, $i=1,\dots,r$ the curve $C$ has the parametrization
of the form $$ u=\tau(a_{0i}+a_{1i}\tau+O(\tau)),\quad
v=v_{1i}+\tau^p (b_{0i}+b_{1i}\tau+O(\tau)), \eqno (21) $$ where
$\tau \in \Complex$ is a local parameter, $a, b\in \Complex$,
$a_{0i}\neq 0$ and $p$ is a positive integer.

$L_2=\prod\limits_{i=1}^k(v-v_{2i})^{m_i}$, $k\leq n$ contains factors of multiplicity $m_i>1$ such
that the corresponding points $(0,v_{2i})$ satisfy the condition
$\displaystyle\frac{\partial F}{\partial u}(0,v_{2i})\neq 0$. For arbitrary $1\leq i \leq k$ we
can write the parametrization of $C$ near $(0,v_{2i})$ as follows
$$
u=\tau^{m_i}(c_{0i}+c_{1i}\tau+O(\tau)),\quad v=v_{2i}+\tau (e_{0i}+e_{1i}\tau+O(\tau))
,\quad c_{0i}, e_{0i} \neq 0. \eqno (22)
$$
At last, the factor $L_3=\prod\limits_{i=1}^s(v-v_{3i})^{l_i}$, $s\leq n$ includes the multipliers
of (20) for which $l_i>1$ and $\displaystyle\frac{\partial F}{\partial u}(0,v_{3i})=0$.

These points are singular and according to Theorem 4 near the point $(v_{3i},0)$ we have
$p_i>1$ local components of $C$ each of them can be parametrized as
$$
u=\tau^{k_{ij}} (g_{0ij}+g_{1ij}\tau+O(\tau)), \quad v=v_{3i}+\tau^{d_{ij}}
,\quad g_{0ij} \neq 0,\quad j=1,\dots,p_i.\eqno (23)
$$
where $d_{ij}$, $k_{ij}$ are positive integers and $\sum\limits_{j=1}^{p_i} k_{ij}\leq l_i$.

In addition we have $$ r+\sum_{i=1}^km_i+\sum_{i=1}^s l_i=n \quad
\mathrm{and} \quad C\cap L_\infty=V_1 \cup V_2 \cup V_3,\quad
V_i=\{L_i=0\},\quad i=1,2,3. $$

From (15) with use of (21), (22), (23) one can show that the following estimates hold
$$
\begin{array}{ll}
\mathrm{deg}(\omega)\mid_{V_1}\leq r(m-2),\quad
\mathrm{deg}(\omega)\mid_{V_2}\leq(m-1)\sum\limits_{i=1}^km_i-k,\\
\mathrm{deg}(\omega)\mid_{V_3}\leq(m-1)\sum\limits_{i=1}^sl_i-\sum\limits_{X\in
 \mathrm{Sing}(C) \cap L_\infty}\mid\pi^{-1}(X)\mid.
\end{array}
$$ Summing we obtain $$ \mathrm{deg}(\omega)\mid_{C\cap
L_\infty}\leq n(m-1)-\sum_{X\in \mathrm{Sing}(C)\cap
L_\infty}\mid\pi^{-1}(X)\mid-k-r. $$ Since
$\mathrm{deg}(\omega)=\mathrm{deg}(\omega)\mid_{C\cap
L_\infty}(\omega)+\mathrm{deg}(\omega)\mid_{C\cap \Complex^2
}(\omega)$ in view of (16), (19) we have

\vspace{0.3cm}

\noindent {\bf Theorem 5.} For an arbitrary algebraic invariant
curve of the system (1) the following estimate for the genus $g$
holds $$ 2g-2\leq n(m-1)-\sum_{X\in
\mathrm{Sing}(C)}\mid\pi^{-1}(X)\mid. \eqno (24) $$

This result seems to be a consequence of the formula 1 of the
paper [2].

\begin{center}
{\bf 5. The algebraic invariant curves with nodes.}
\end{center}

Let $C$ be an algebraic invariant curve of the system (1) with the defining polynomial $f(x,y)$.
\vspace{0.3cm}

\noindent{\bf Lemma 4.} $\mid \mathrm{Sing}(C)\mid\leq
m^2+\displaystyle\frac{n}{2}$. \vspace{0.3cm}

This is a simple consequence of Corollary 2 and the notation, that $C$ has at most $n/2$ singular
points at infinity.
\vspace{0.3cm}

\noindent {\bf Theorem 6.} Let there exists the integer $K$ such
that $\forall$ $ X \in \mathrm{Sing}(C)$ we have
$\left(f,\frac{\partial f}{\partial y}\right)_X\leq K$, then the
following estimate for the degree of the curve $C$ holds $$
n\leq\displaystyle\frac{4+2m+K+((4+2m+K)^2+16Km^2)^{1/2}}{4},
\eqno (25) $$ where $m$ is the degree of the system (1).
\vspace{0.3cm}

\noindent{\bf Proof.} With using of (17), (24) one can show that
$$ n(n-3)-\sum_{X\in \mathrm{Sing}(C)}\left(f,\frac{\partial
f}{\partial y}\right)_X \leq n(m-1).\eqno (26) $$ By our
assumption: $\left(f,\frac{\partial f}{\partial y}\right)_X\leq
K$. According to Lemma 4 we obtain immediately $$ \sum_{X\in
\mathrm{Sing}(C)}\left(f,\frac{\partial f}{\partial
y}\right)_X\leq K(m^2+\displaystyle\frac{n}{2}). \eqno (27) $$
Puting (27) into (26) we arrive at Theorem 6. \vspace{0.3cm}

\noindent{\bf Corollary 3.} Let us suppose that all singular
points of the algebraic invariant curve $C$ are nodes, then $$
n\leq 2(m+1). \eqno (28) $$ \vspace{0.3cm} Indeed, as a node is an
ordinary double point then $K=1$ and we can use the estimate (25)
which gives (28). It is interesting to compare this result with
Theorem 3 of the paper [2].

\begin{center}
{\bf Acknowledgement.}
\end{center}

The author would be grateful to L. Gavrilov  for his attention to
the paper and many useful comments.

\begin{center}

{\bf References}

\end{center}

 \noindent [1] A. Campillo, M.M. Carnicer, {\it Proximity
inequalities and bounds for the degree of invariant curves by
foliations of $P^2_C$}, Trans. Amer. Math. Soc. 349, (1997), no.
6, 2211--2228

\noindent [2] D. Cerveau, A. Lins Neto, {\it Holomorphic
foliations in $CP(2)$ having an invariant algebraic curve}, Ann.
Inst. Fourier, Grenoble, 41, 4 (1991), 883--903

\noindent [3] Chavarriga, Javier; Llibre, Jaume, {\it On the
algebraic limit cycles of quadratic systems}. Proceedings of the
IV Catolan Days of Applied Mathematics (Tarragona, 1998), 17-24,
Univ. Rovira virgili, Tarragona, 1998.

\noindent [4] G. Darboux,{ \it M\'emoire sur les \'equations
diff\'erentielles alg\'ebrique du premier ordre et du premier
degr\'e (M\'elanges)}, Bull. Sci. Math. (1878), 60-96; 123-144;
151-200.

\noindent [5] M. Galeotti, {\it Monodromic Unbounded Polycycles},
Annali di Matematica pura ed applicata (IV), Vol. CLXXI (1996),
83-85

\noindent [6] E. Gonzales-Velasco, {\it Generic properties of
polynomial vector fields at infinity}, Trans. Am. Math. Soc., 143
(1968)

\noindent [7] P. A. Griffits, {\it Introduction to algebraic
curves}, Transactions of mathematical monographs 76, American
Mathematical Society (1989)

\noindent [8] D.Hilbert, {\it Mathematical problems}, Bull. Amer.
Math. Soc. 8 (1902), 437-479.

\noindent [9] A. Jablonskii, {\it Algebraic integrals of system of
differential equations}, (Russian, Engl. transl.), Diff. Uravn. 6
no. 11, 1970, 1752-1760; Engl. transl. 1326-1333.

\noindent [10]  R. E. Kooij, C. J. Christopher, {\it Algebraic
invariant curves and the integrability of polynomial systems},
Appl. Math. Lett. 6 (1993), 51-53

\noindent [11] J. Moulin--Ollagnier, A. Nowicki, J.-M. Strelcyn,
{\it On the non-existence of constants of derivations: the proof
of theorem of Jouanolou and its development}, Bull. Sci. math. 119
(1995), 195-233.

\noindent [12] H.Poincar\'e, "Sur les courbes d\'efinies par les
\'equations diff\'erentielles",
 Oeuvres de Henri Poicar\'e, Paris, Gauthiers-Villars et Cie, Editeurs, 1928, vol. 1

\vspace{0.3cm}

 \noindent{\small \bf
 Section de Mathematiques,\\
 Universit\'e de Gen\`eve\\
 2-4, rue du Lievre,\\
 CH-1211, Case postale 240, Suisse \\
Tel l.: +41 22 309 14 03 \\
 Fax: +41 22 309 14 09 \\
E--mail: Alexei.Tsygvintsev@math.unige.ch}

\end{document}